\documentclass[pdflatex]{article}
\usepackage{graphicx} 
\usepackage{amsmath}
\usepackage{amssymb}
\usepackage{amsfonts}
\usepackage{todonotes}
\usepackage{subcaption}
\usepackage{hyperref}
\usepackage{lineno}
\usepackage{setspace}
\usepackage{geometry}
\usepackage{nomencl}
\makenomenclature

\usepackage{natbib}

\title{Discrete Exterior Calculus-based Whittle-Matérn Fields on Meshes}
\author{Michael Gillan\footnote{University of Exeter -- m.gillan@exeter.ac.uk}}
\date{}

\begin{document}

\maketitle

\begin{abstract}
    
This work presents a sparse precision formulation for Whittle-Matérn fields based on Discrete Exterior Calculus. By treating the $k$-cochain as the latent object to be modeled, we show that a single formulation can produce Whittle-Matérn fields for all $k$ and integer $\alpha$ values. By utilising Whitney maps we can enable a push-forward of the prior on 1-cochains, enabling inference for discrete vector fields. 
We demonstrate this model on interpolation of vector fields on a sphere.
\end{abstract}

\section{Introduction}

Gaussian Processes (GPs) are a commonly-used tool for spatial modelling under uncertainty. A GP is a multivariate distribution where the joint distribution of every finite set of variables is Gaussian \citep{rasmussenGaussianProcessesMachine2006}. A GP is wholly defined by its mean and covariance functions, and the mean and variance of an unobserved spatial index is defined by the GP and any observed data. An \emph{isotropic} GP is one where the covariance function is a function of the Euclidean distance between two spatial locations, and common GPs, such as the Matérn or squared exponential GPs, are isotropic.

However, on curved manifolds like spheres, Euclidean distance functions are not sufficient for defining spatial covariance, and using geodesic distance instead generally does not produce positive definite covariance. \cite{borovitskiyMaternGaussianProcesses2020} introduce an alternative formulation for scalar covariance kernels defined on an arbitrary Riemannian manifold $(M, g)$ through the spectral decomposition of the Laplace-Beltrami operator $-\Delta_g$. Both \cite{hutchinsonVectorvaluedGaussianProcesses2021a} and \cite{robert-nicoudIntrinsicGaussianVector2024b} extend this idea to modelling vector fields, the former through projected kernels, and the latter through applying derivative operators to the eigenfunctions of the Laplace-Beltrami to get vector fields. In particular, \cite{robert-nicoudIntrinsicGaussianVector2024b} show that \emph{Hodge-compositional} processes can be constructed in this way, allowing separate parametrisations of the curling and diverging velocities in a vector field. \cite{borovitskiyMaternGaussianProcesses2021b} show that the spectral formulation of the Matérn GP can be applied to graphs and triangular meshes through the eigenpairs of the graph Laplacian, and show empirically that it converges to its smooth counterpart in the Riemannian setting. \cite{yangHodgeCompositionalEdgeGaussian2024a} propose a methodology for modelling functions on the edges of meshes, extending the formulation of intrinsic vector GPs with the graph Laplacian. \cite{gillanDiscreteGaussianVector2026} show that the Hodge-compositional vector field model of \cite{robert-nicoudIntrinsicGaussianVector2024b} has a discrete analogue by substituting the discrete exterior calculus (DEC) cotangent Laplacian for $-\Delta_g$ and showing that it can be used to model vector fields defined on the vertices of a manifold mesh. However, this precludes a sparse formulation. \cite{alainGaussianProcessesCellular2024} develop a model for Matérn GPs over the components of a cellular complex, analogous to the simplicial complex of DEC, and propose a kernel for propagating information between different orders of cells. \cite{alain2026probabilistic} develop a similar model to \cite{gillanDiscreteGaussianVector2026} using \emph{finite element exterior calculus}, eliciting a similar covariance-based Hodge-compositional model without the restriction to mesh vertices.

Standard dense Gaussian-process calculations require O($N^3$) factorisation and O($N^2$) storage. The SPDE approach \citep{lindgrenExplicitLinkGaussian2011} explicitly constructs Matérn precision matrices for Gaussian Markov Random Fields (GMRFs) with respect to a mesh by discretising the SPDE using piece-wise linear basis functions, enabling sparse linear-algebra methods. As highlighted in a later review \citep{lindgrenSPDEApproachGaussian2022}, this result can be generalised to smooth manifolds by replacing the Euclidean $\Delta$ with the equivalent manifold operator $\Delta_g$\nomenclature{$\Delta_g$}{Laplacian defined on a Riemannian Manifold $M$ with metric $g$}. \cite{nambirajanApproximatingGaussianWhittleMatern2026} show that a Gaussian Markov Random Field (GMRF) approximation for Whittle-Matérn fields can be built from DEC operators for fields on vertices and faces.

\paragraph{Contributions} This work presents a sparse precision formulation based on DEC for Whittle-Matérn fields on $k$-forms. We show that, through operators defined with respect to a simplicial mesh of a manifold, a precision matrix can be constructed over $k$-cochains for every $k$ and integer $\alpha$, generalising previous constructions for $k=0$ and $k=n$. The prior and posterior precision matrices are sparse, allowing the use of sparse linear-algebra methods for posterior inference. This model can be applied directly to point, edge, and areal data, with observations and prediction locations anywhere on the mesh through Whitney maps, including vector field data through a Whitney observation operator. We show how to use this formulation to perform posterior inference, and demonstrate it for vector-field inference on the sphere.


\section{Background}

\subsection{Discrete Exterior Calculus}

Exterior calculus is the calculus defined by generalising the derivative to operate on higher-degree differential forms \citep{flemingFunctionsSeveralVariables1987}. A differential $k$-form is a quantity that can be integrated over a $k$-dimensional manifold; for example, the expression $h(x) dx$ is a 1-form and can be integrated over a 1-dimensional line segment by $\int_a^b h(x) dx$. The discretisation of exterior calculus is called discrete exterior calculus (DEC) \citep{hiraniDiscreteExteriorCalculus2003a}. An implementation of most of the core theory is provided in the PyDEC library \citep{bellPyDECSoftwareAlgorithms2012}. 

In the following, $\omega, \eta$ are smooth $k$-forms, $\alpha, \beta$ are discretised $k$-forms, cochains, $c$ is a $k$-chain, $M$ is the Riemannian manifold that is triangulated by the simplicial complex $K$, and $\sigma^k$ is a simplex of dimension $k$. For an in-depth treatment see, e.g. \cite{desbrunDiscreteExteriorCalculus2005}.

\subsubsection{Simplicial Complexes, Chains and Cochains}

A smooth Riemannian manifold $(M, g)$\nomenclature{$M$}{Riemannian manifold with metric $g$} of dimension $n$ can be discretised as a finite collection of 0-, 1-, \dots, $(n-1)$-, and $n$-dimensional piecewise linear components. For a 2-dimensional manifold, these would be vertices (0-dimensional), edges (1-dimensional), and faces (2-dimensional). This collection of components is called a simplicial complex, $K$\nomenclature{$K$}{Simplicial complex}, and the individual $k$-dimensional components, $\sigma^k$\nomenclature{$\sigma^k$}{A $k$-dimensional simplex}, are called $k$-simplices. We assume $K$ triangulates $M$ and that $K$ is embedded in Euclidean space with vertex positions chosen so that its lengths, areas and angles approximate those of $(M,g)$. We write $K^{(k)}$ for the set of $k$-simplices (e.g. $K^{(1)}$ is the collection of edges on a triangulated surface) and $N_k=|K^{(k)}|$.

The orientation of a $k$-simplex is the ordering of its vertices $[v_0,\dots,v_k]$, up to even permutations. The \emph{chain group} $C_k(K)$ is the abelian group generated by the oriented $k$-simplices, where opposite orientations are negative. Fixing one orientation per simplex, $C_k(K)$ is the free abelian group on $K^{(k)}$,
\begin{equation}
  C_k(K)=\Big\{\textstyle\sum_{\sigma^k\in K^{(k)}}c_{\sigma^k}\,\sigma^k : c_{\sigma^k}\in\mathbb{Z}\Big\}
  \cong\mathbb{Z}^{N_k},
\end{equation}
and we set $C_k(K)=0$ for $k<0$ or $k>n$. The orientations of individual simplices are an arbitrary choice.

A \emph{discrete $k$-form}, or $k$-cochain, is a homomorphism $\alpha:C_k(K)\to\mathbb{R}$, and these form the \emph{cochain group}
\begin{equation}
  C^k(K):=\operatorname{Hom}(C_k(K),\mathbb{R}).
\end{equation}
A cochain is determined by its values $\alpha(\sigma^k)$ on the simplices, and $C^k(K)\cong\mathbb{R}^{N_k}$; a $0$-cochain is a function on vertices, a $1$-cochain a value on each oriented edge, and so on. We write
\begin{equation}
  \langle\alpha,c\rangle:=\alpha(c)=\sum_{\sigma^k\in K^{(k)}}\alpha(\sigma^k)\,c_{\sigma^k},
  \qquad \alpha\in C^k(K),\; c\in C_k(K),
\end{equation}
for the evaluation of a cochain on a chain. 

The circumcentric dual of $K$ associates each $k$-simplex $\sigma^k$ with an $(n-k)$-dimensional dual cell $\star\sigma^k$ built from the circumcentres of the simplices containing $\sigma^k$. On a triangulated surface, the dual of a face is its circumcentre, the dual of an edge is the segment joining the circumcentres of its adjacent faces, and the dual of a vertex is the polygon joining the circumcentres of its incident faces. We write $|\sigma^k|$ and $|\star\sigma^k|$ for the $k$- and $(n-k)$-dimensional volumes of a simplex and its dual cell, taking the volume of a point to be 1. On well-centred meshes (every simplex contains its circumcentre) all dual volumes are positive.

\subsubsection{de Rham and Whitney Maps}

Smooth and discrete forms are related through the \emph{de Rham} and \emph{Whitney} maps \citep{whitneyGeometricIntegrationTheory1957}. The de Rham map $\psi^k$\nomenclature{$\psi^k$}{de Rham map} discretises a smooth $k$-form $\omega$ by integration over $k$-simplices,
\begin{equation}
    \psi^k(\omega)(\sigma^k) = \langle\psi^k(\omega), \sigma^k\rangle := \int_{\sigma^k} \omega.
\end{equation}
$\psi^k\omega$ is a well-defined cochain in $C^k(K)$, and by linearity $\langle\psi^k\omega,c\rangle=\int_c\omega$ for every chain $c\in C_k(K)$.

In the opposite direction, the Whitney map $\phi^k$\nomenclature{$\phi^k$}{Whitney map} interpolates a cochain $\alpha\in C^k(K)$ as a piecewise-smooth $k$-form,
\begin{equation}\label{eq:whitney-map}
  \phi^k\alpha = \sum_{\sigma^k\in K^{(k)}}\alpha(\sigma^k)\,W_{\sigma^k} .
\end{equation}
Writing $\lambda_i$ for the barycentric coordinate of vertex $v_i$, equal to one at $v_i$ and zero at all other vertices, the Whitney forms of degree zero, one and two are
\begin{equation}\label{eq:whitney-forms}
  W_{v_i}=\lambda_i,
  \quad
  W_{[v_i,v_j]}=\lambda_i\,\mathrm{d}\lambda_j-\lambda_j\,\mathrm{d}\lambda_i,
  \quad
  W_{f} = \frac{\mathbf{1}_f}{|f|} \mathrm{d}A_f,
\end{equation}
where $f = [v_i, v_j, v_l]$ is a face on a triangulated surface and $\mathrm{d}A_f$ is its area form, oriented by the vertex ordering of $f$. These satisfy $\int_{\sigma^k} W_{\sigma'^k}=\delta_{\sigma^k\sigma'^k}$ for all $\sigma,\sigma'\in K^{(k)}$, so that
\begin{equation}
  \psi^k\phi^k = I,
\end{equation}
i.e. interpolating a cochain and then discretising recovers it exactly.

\subsubsection{Exterior Derivative}

The discrete boundary operator $\partial_k:C_k(K)\to C_{k-1}(K)$ is defined by
\begin{equation}
  \partial_k[v_0,\dots,v_k]=\sum_{i=0}^{k}(-1)^i[v_0,\dots,\hat v_i,\dots,v_k],
\end{equation}
where $[v_0, \dots, \hat v_i, \dots, v_k]$ is the $(k-1)$-dimensional simplex obtained by omitting the vertex $v_i$ \citep{desbrunDiscreteExteriorCalculus2005}. $\partial_k$ is represented by a sparse signed incidence matrix with
entries in $\{0,\pm1\}$, and $\partial_k\partial_{k+1} = 0$.

The \emph{discrete exterior derivative}, $\mathbf{d}_k : C^k(K) \to C^{k+1}(K)$, \nomenclature{$\mathbf{d}_k$}{Discrete exterior derivative on $k$-cochains} is the discrete difference operator acting on $k$-cochains. For any $(k+1)$-chain $c$, Stokes' theorem can be written as an adjoint relationship 

\begin{equation}
    \int_c \mathrm{d} \omega = \int_{\partial c}\omega \iff \langle  \psi^{k+1}\mathrm{d}\omega, c \rangle = \langle \psi^k\omega, \partial_{k+1} c \rangle.
\end{equation}
We define $\mathbf{d}$ as the adjoint of the discrete boundary operator,

\begin{equation}
    \langle \mathbf{d}_k \alpha, c \rangle = \langle \alpha, \partial_{k+1} c \rangle,
  \qquad \alpha\in C^k(K),\; c\in C_{k+1}(K).
\end{equation}
It follows from $\partial_k\partial_{k+1} = 0$ that $\mathbf{d}_{k+1}\mathbf{d}_k = 0$. Additionally, we can see that $\mathbf{d}_k$ commutes with the de Rham map $\psi^k$,
\begin{equation}
    \langle \mathbf{d}_k \psi^k \omega, c \rangle = \langle \psi^k \omega, \partial_{k+1} c \rangle = \int_{\partial c} \omega = \int_c \mathrm{d}\omega = \langle \psi^{k+1} \mathrm{d}\omega, c \rangle,
\end{equation}
for any smooth $k$-form $\omega$ and $c\in C_{k+1}(K)$. Since a cochain is determined by its values on chains, $\mathbf{d}_k\psi^k=\psi^{k+1}\mathrm{d}$, and $\mathbf{d}_k$ introduces no discretisation error. In matrix form, $\mathbf{d}_k$ is the transpose of the signed incidence matrix of $\partial_{k+1}$, and is simply a sparse difference matrix as a result of the sparsity of $\partial_{k+1}$.

\subsubsection{Hodge Star}

In the DEC framing, all metric information enters through an inner product on cochains. For smooth $k$-forms on $(M, g)$, the $L^2$ inner product is
\begin{equation}
  (\omega,\eta)_{L^2} = \int_M g(\omega,\eta)\,\mathrm{d}V,
  \qquad
  g(\omega,\eta) = \sum_{i_1<\dots<i_k}
  \omega(e_{i_1},\dots,e_{i_k})\,\eta(e_{i_1},\dots,e_{i_k}),
\end{equation}
where $\{e_i\}$ is any local orthonormal frame. We seek a symmetric positive definite matrix $\star_k$ such that
\begin{equation}
    (\psi^k\omega)^\top\star_k(\psi^k\eta)\approx(\omega,\eta)_{L^2},
\end{equation}
and define the
discrete inner product
\begin{equation}\label{eq:hodge-ip}
  (\alpha,\beta)_k := \alpha^\top\star_k\,\beta.
\end{equation}
The natural choice is the Gram matrix of the Whitney forms, $(\star^W_k)_{\sigma\sigma'}=(W_\sigma,W_{\sigma'})_{L^2}$, for which $(\alpha,\beta)_k=(\phi^k\alpha,\phi^k\beta)_{L^2}$ exactly; $\star^W_k$ is then the finite-element mass matrix of $k$-cochains. 

In practice, DEC uses a diagonal (lumped) approximation of this mass matrix. For $\alpha = \psi^k\omega$, $\alpha(\sigma^k)=\int_{\sigma^k}\omega\approx\bar\omega_{\sigma^k}|\sigma^k|$, and the ratio $\alpha(\sigma^k)/|\sigma^k|$ approximates the average density of $\omega$ along $\sigma^k$. Giving $\sigma^k$ a region of measure proportional to $|\sigma^k|\,|{\star\sigma^k}|$ and summing density squared times volume gives 
\begin{equation}
  (\omega,\omega)_{L^2}
  \approx\sum_{\sigma^k\in K^{(k)}}
  \left(\frac{\alpha(\sigma^k)}{|\sigma^k|}\right)^2|\sigma^k|\,|{\star\sigma^k}|
  =\sum_{\sigma^k\in K^{(k)}}\frac{|{\star\sigma^k}|}{|\sigma^k|}\,\alpha(\sigma^k)^2 = \alpha^\top \star_k \alpha, 
\end{equation}
where $\star_k$ is the diagonal Hodge star,
\begin{equation}
  \star_k=\operatorname{diag}\!\left(\frac{|{\star\sigma^k}|}{|\sigma^k|}\right)_{\sigma^k\in K^{(k)}}.
\end{equation}
$\star_k$ is symmetric positive definite on well-centred meshes, where all $|\star \sigma^k|$ are positive.

\subsubsection{Co-derivative}

The \emph{discrete co-derivative}
$\delta_{k+1}:C^{k+1}(K)\to C^k(K)$ is defined as the adjoint of
$\mathbf{d}_k$ with respect to the inner product in Eq. \ref{eq:hodge-ip},
\begin{equation}
  (\mathbf{d}_k\alpha,\beta)_{k+1}=(\alpha,\delta_{k+1}\beta)_k
  \qquad\text{for all }\alpha\in C^k(K),\;\beta\in C^{k+1}(K).
\end{equation}
Since $\star_k$ is symmetric positive definite, and hence invertible,
\begin{equation}
  (\mathbf{d}_k\alpha,\beta)_{k+1}
  =\alpha^\top\mathbf{d}_k^\top\star_{k+1}\beta
  =\alpha^\top\star_k\left(\star_k^{-1}\mathbf{d}_k^\top\star_{k+1}\beta\right),
\end{equation}
so
\begin{equation}
  \delta_{k+1}=\star_k^{-1}\mathbf{d}_k^\top\star_{k+1}.
\end{equation}
Just as $\mathbf{d}_k$ is the adjoint of $\partial_{k+1}$ under the pairing $\langle\cdot,\cdot\rangle$, $\delta_{k+1}$ is the adjoint of $\mathbf{d}_k$ under the inner product $(\cdot,\cdot)_k$.

\subsubsection{Hodge-Laplacian}

The \emph{discrete Hodge-Laplacian}\footnote{also called the Laplace-de Rham operator in some literature} on $k$-cochains is

\begin{equation} \label{eq:hodge-lap}
    \Delta_k = \mathbf{d}_{k-1} \delta_k + \delta_{k+1}\mathbf{d}_k,
\end{equation}
\nomenclature{$\Delta_k$}{Discrete Hodge-Laplacian\nomrefeq}
with $\delta_0=0$ and $\mathbf{d}_n=0$, since $C^{-1}(K)=C^{n+1}(K)=0$. $\Delta_k$ is self-adjoint and positive semidefinite with respect to $(\cdot,\cdot)_k$, since $(\Delta_k\alpha,\alpha)_k=(\mathbf{d}_k\alpha,\mathbf{d}_k\alpha)_{k+1} +(\delta_k\alpha,\delta_k\alpha)_{k-1}\geq0$. $\Delta_k$ as a matrix is not symmetric, but $\star_k\Delta_k$ is.

\subsection{Whittle-Matérn Fields} 


The Matérn covariance function for locations $\mathbf{u}, \mathbf{v} \in \mathbb{R}^n$ and Euclidean distance $\|\cdot\|$ is defined as
\begin{equation}
    r(\mathbf{u, v}) = \frac{\sigma^2}{2^{\nu - 1} \Gamma(\nu)}(\kappa \|\mathbf{v - u}\|)^\nu K_\nu (\kappa \|\mathbf{v - u}\|),
\end{equation}
where $K_\nu$ is the modified Bessel function of the second kind with order $\nu > 0$, $\kappa > 0$ is a scaling parameter and $\sigma^2$ is the marginal variance. The Matérn random field has mean-square derivatives of integer orders $m < \nu$. \cite{lindgrenExplicitLinkGaussian2011} suggests a more natural interpretation of $\kappa$ as a range parameter $\rho$, with an empirical reparametrisation $\rho = \frac{\sqrt{8\nu}}{\kappa}$ giving correlations near 0.1 when $\|\mathbf{v - u}\| = \rho$.

\cite{whittleStochasticProcessesSeveral1963} formalises the Matérn GP as the solution to the SPDE 
\begin{equation}\label{eq:WM-SPDE}
(\kappa^2 - \Delta)^{\alpha/2}x(\mathbf{s}) = \mathcal{W}(\mathbf{s}),    
\end{equation}
with the same $\kappa > 0$ and $\alpha = \nu + \frac{n}{2}$, and a white noise process $\mathcal{W}$\nomenclature{$\mathcal{W}$}{White noise process}. $\Delta$ is the Euclidean Laplacian\nomenclature{$\Delta$}{Euclidean Laplacian} and the marginal variance of the field is given by 
\begin{equation}
    \sigma^2 = \frac{\Gamma(\nu)}{\Gamma(\nu + n/2)(4\pi)^{n/2}\kappa^{2\nu}}.    
\end{equation}

\section{Methodology}

\subsection{Whittle-Matérn Fields on $k$-Cochains} \label{sec:wmf}

When a Gaussian vector is defined implicitly by a linear system driven by Gaussian noise, its precision is defined by the system matrix and the noise covariance. Suppose
\begin{equation}\label{eq:linear-gaussian}
  Sx=b,\qquad b\sim\mathcal{N}(0,\Sigma_b),
\end{equation}
with $S$ invertible. Then $x=S^{-1}b$ is Gaussian with
\begin{equation}
  \operatorname{Cov}(x)=S^{-1}\Sigma_b S^{-\top},
  \qquad
  Q=\operatorname{Cov}(x)^{-1}=S^\top\Sigma_b^{-1}S.
\end{equation}

We can translate the Whittle-Matérn SPDE formulation in Eq. \ref{eq:WM-SPDE} to DEC by replacing $-\Delta$ with the $k$-form Hodge-Laplacian $\Delta_k$ as defined in Eq. \ref{eq:hodge-lap}.
\begin{equation} \label{eq:DEC-SPDE}
    (\kappa^2 + \Delta_k)^\frac{\alpha}{2} x = \mathcal{W}.
\end{equation}
Note the change of sign between Eqs. \ref{eq:WM-SPDE} and \ref{eq:DEC-SPDE} since $\Delta$ is negative semi-definite and $\Delta_k$ is positive semi-definite.

To discretise Eq. \ref{eq:DEC-SPDE} for $\alpha=2$, we use its weak form: we require $\big((\kappa^2+\Delta_k)x,v\big)_k=(\mathcal{W},\phi^kv)_{L^2}$ for all $v\in C^k(K)$. For each $\sigma\in K^{(k)}$, the left-hand side is the $\sigma$-th entry of $\star_k(\kappa^2+\Delta_k)x=Sx$, and the right-hand side is $b_{\sigma}=(\mathcal{W},W_\sigma)_{L^2}$, giving $Sx=b$. Since
$\kappa^2+\Delta_k$ is self-adjoint with respect to $(\cdot,\cdot)_k$, $S$ is symmetric, and its explicit form follows from the quadratic form $u^\top Su=((\kappa^2+\Delta_k)u,u)_k$,
\begin{align}
\begin{split}
    ((\kappa^2 + \Delta_k) u, u)_k& = ((\kappa^2 + \mathbf{d}_{k-1} \delta_k + \delta_{k+1}\mathbf{d}_k) u, u )_k \\
    & = \kappa^2 (u, u)_k + (\mathbf{d}_{k-1} \delta_k u, u)_k + (\delta_{k+1}\mathbf{d}_k u, u )_k \\
    & = \kappa^2 (u, u)_k + (\delta_k u, \delta_k u )_{k-1} + (\mathbf{d}_k u, \mathbf{d}_k u)_{k+1} \\
    & = \kappa^2 \|u\|^2_k + \|\delta_k u\|^2_{k-1} + \|\mathbf{d}_k u\|_{k+1}^2\\
    & = \kappa^2 u^\top \star_k u + (\delta_k u)^\top \star_{k-1} (\delta_k u) + (\mathbf{d}_k u)^\top \star_{k+1} (\mathbf{d}_k u) \\
    & = \kappa^2 u^\top \star_k u + u^\top \star_k \mathbf{d}_{k-1} \star^{-1}_{k-1} \mathbf{d}^\top_{k-1} \star_k u + u^\top \mathbf{d}^\top_{k} \star_{k+1} \mathbf{d}_k u \\
    & = u^\top S u,\\
    &\text{where} \\
    S &= \kappa^2 \star_k + \star_k \mathbf{d}_{k-1} \star^{-1}_{k-1} \mathbf{d}^\top_{k-1} \star_k + \mathbf{d}^\top_k \star_{k+1}\mathbf{d}_k \\
    &= \star_k (\kappa^2 + \Delta_k).
\end{split}
\end{align}
$S\in\mathbb{R}^{N_k\times N_k}$ is sparse since $\star_{k-1}^{-1}$ is diagonal. Using $\star_k^W$ would make the term $\star_k \mathbf{d}_{k-1} \star^{-1}_{k-1} \mathbf{d}^\top_{k-1} \star_k$ dense, which justifies our use of the lumped mass matrix.

The white noise process $\mathcal{W}$ is a random linear functional such that
\begin{equation}
    \mathbb{E}[(\mathcal{W}, v)_{L^2} (\mathcal{W}, w)_{L^2}] = (v, w)_{L^2}.
\end{equation}
Using Whitney forms $\{W_{\sigma}\}$ as a basis, $b_\sigma := (\mathcal{W}, W_\sigma)_{L^2}$, and
\begin{equation}
    \text{Cov}(b)_{\sigma\sigma'} = \mathbb{E}[b_\sigma b_{\sigma'}] = \mathbb{E}[(\mathcal{W}, W_\sigma)_{L^2} (\mathcal{W}, W_{\sigma'})_{L^2} ] = (W_\sigma, W_{\sigma'})_{L^2} = \int_\mathcal{M} g(W_\sigma, W_{\sigma'}) dV = (\star^W_k)_{\sigma\sigma'}.
\end{equation}
Approximating $\star^W_k$ by the lumped mass matrix $\star_k$ gives $b \sim \mathcal{N}(0, \star_k)$. Thus, when $\alpha = 2$, the weak form of Eq. \ref{eq:DEC-SPDE} is 
\begin{equation}
S x = b, \quad b \sim \mathcal{N}(0, \star_k),     
\end{equation}
and the precision matrix of $x$ is given by
\begin{equation}\label{eq:q_2}
    Q_2 = S \star^{-1}_k S.
\end{equation}

For general $\alpha \neq 2$, the fractional power in Eq. \ref{eq:DEC-SPDE} is most naturally defined after \emph{whitening}. The discrete inner product satisfies
\begin{equation}
  (u,v)_k = u^\top\star_k v = (\star_k^{1/2}u)^\top(\star_k^{1/2}v),
\end{equation}
so the change of variables $\tilde x=\star_k^{1/2}x$ transforms the $\star_k$ inner product to the Euclidean one. The self-adjoint operator $\kappa^2+\Delta_k$ is then represented by the symmetric positive definite matrix
\begin{equation}
  \tilde S=\star_k^{1/2}(\kappa^2+\Delta_k)\star_k^{-1/2}
          =\star_k^{-1/2}S\,\star_k^{-1/2},
\end{equation}
using $S=\star_k(\kappa^2+\Delta_k)$. The noise process becomes
\begin{equation}
    \tilde{b} = \star_k^{-\frac{1}{2}} b \sim \mathcal{N}(0, \star_k^{-\frac{1}{2}} \star_k \star_k^{-\frac{1}{2}}) = \mathcal{N}(0, I).
\end{equation}

For $\alpha=2$, $Sx=b$ is equivalent to $\tilde S\tilde x=\tilde b$. For general $\alpha$ the whitened form of Eq. \ref{eq:DEC-SPDE} is
\begin{equation}
  \tilde S^{\alpha/2}\tilde x=z,\qquad z\sim\mathcal{N}(0,I).
\end{equation}
Hence $\operatorname{Cov}(\tilde x)=\tilde S^{-\alpha}$. Transforming back
via $x=\star_k^{-1/2}\tilde x$ gives
$\operatorname{Cov}(x)=\star_k^{-1/2}\tilde S^{-\alpha}\star_k^{-1/2}$, and
so
\begin{equation}
  Q_\alpha=\star_k^{1/2}\,\tilde S^{\alpha}\,\star_k^{1/2}.
\end{equation}

For integer $\alpha$, we can use telescoping to define $Q_\alpha$,

\begin{align} \label{eq:q_alpha}
\begin{split}        
    \tilde{S}^\alpha &= \underbrace{(\star_k^{-\frac{1}{2}} S \star_k^{-\frac{1}{2}})(\star_k^{-\frac{1}{2}} S \star_k^{-\frac{1}{2}}) \dots (\star_k^{-\frac{1}{2}} S \star_k^{-\frac{1}{2}})}_{\alpha \text{ copies}} \\
    &= \star_k^{-\frac{1}{2}} \underbrace{S \star_k^{-1} S \star_k^{-1} \dots S}_{\alpha \text{ copies of } S} \star_k^{-\frac{1}{2}} \\
    Q_\alpha &= \star_k^\frac{1}{2} \star_k^{-\frac{1}{2}} \underbrace{S \star_k^{-1} S \star_k^{-1} \dots S}_{\alpha \text{ copies of }S} \star_k^{-\frac{1}{2}} \star_k^\frac{1}{2} \\
    &= \underbrace{S \star_k^{-1} S \star_k^{-1} \dots S}_{\alpha \text{ copies of } S} \\
    &= S(\star_k^{-1} S)^{\alpha - 1},\quad\alpha\in\mathbb{N}
\end{split}
\end{align}

For $\alpha=1$ this is $Q_1=S$, and for $\alpha=2$ it recovers $Q_2=S\star_k^{-1}S$. For non-integer $\alpha$, $\tilde S^\alpha$ remains well defined through spectral powers but is dense, so $Q_\alpha$ loses both its sparsity and the Markov property. We therefore restrict attention to integer $\alpha$; fractional values can be handled by rational approximation \citep{bolinRationalSPDEApproach2020}.

The smoothness parameter is related to $\alpha$ by $\nu=\alpha-n/2$, and $\nu>0$ requires $\alpha>n/2$. With integer $\alpha$, $\nu$ takes integer values when $n$ is even and half-integer when $n$ is odd. On a surface ($n=2$), $\alpha=1$ gives $\nu=0$, and so does not define a Matérn covariance function, but is still a well-defined random measure \citep{lindgrenExplicitLinkGaussian2011}.

\section{Implementation}

\subsection{Prior Sampling}\label{sec:prior}

To sample from $\mathcal{N}(0, Q_\alpha^{-1})$ we require $L$ such that $LL^\top = Q_\alpha$. We consider the cases where $\alpha = 2m$ and $\alpha = 2m+1$ separately.

When $\alpha = 2m$,
\begin{align}
\begin{split}
Q_\alpha &= S (\star_k^{-1} S)^{2m-1} \\
&=  S (\star_k^{-1} S)^{2m-2} \star_k^{-1} S \\
& = (S (\star_k^{-1} S)^{m-1} \star_k^{-\frac{1}{2}}) (S (\star_k^{-1} S)^{m-1} \star_k^{-\frac{1}{2}})^\top \\
L &= S (\star_k^{-1} S)^{m-1} \star_k^{-\frac{1}{2}}.
\end{split}
\end{align}
When $\alpha = 2m+1$,
\begin{align}
\begin{split}
Q_\alpha &= (S\star_k^{-1})^m S (\star_k^{-1} S)^m \\
& = ((S\star_k^{-1})^{m} L_S)((S\star_k^{-1})^{m} L_S)^\top \\
L & = (S\star_k^{-1})^{m} L_S,
\end{split}
\end{align}
where
\begin{equation}
L_S L_S^\top = S.
\end{equation}
A sample $x$ is then given by $x = L^{-\top} z$, $z \sim \mathcal{N}(0, I)$, i.e. 
\begin{equation}
x = 
\begin{cases}
(S^{-1} \star_k)^{m-1} S^{-1} \star_k^\frac{1}{2} z & \alpha = 2m \\
(S^{-1} \star_k)^m  L_S^{-\top} z & \alpha = 2m+1
\end{cases}
\end{equation}
Figure \ref{fig:samples} shows samples of $k$-cochains on a unit sphere.
\newcommand{\subfigwidth}{0.4}
\begin{figure*}[h!]
        \centering
        \begin{subfigure}[b]{\subfigwidth\textwidth}
            \centering
            \includegraphics[width=\textwidth]{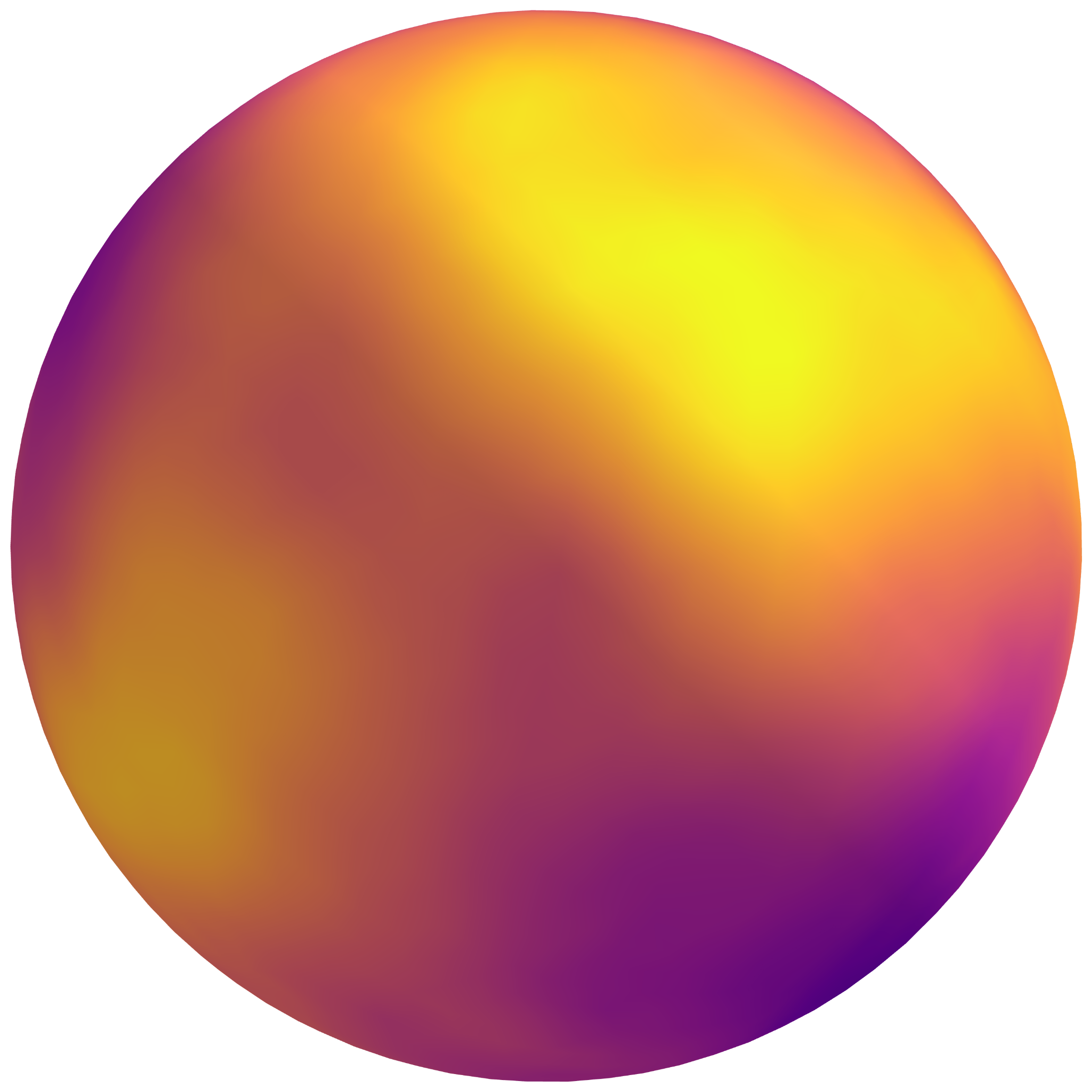}
            \caption{{\small Sample of $0$-cochain (function over points)}}  
            \label{fig:0-sample}
        \end{subfigure}
        \hfill
        \begin{subfigure}[b]{\subfigwidth\textwidth}  
            \centering 
            \includegraphics[width=\textwidth]{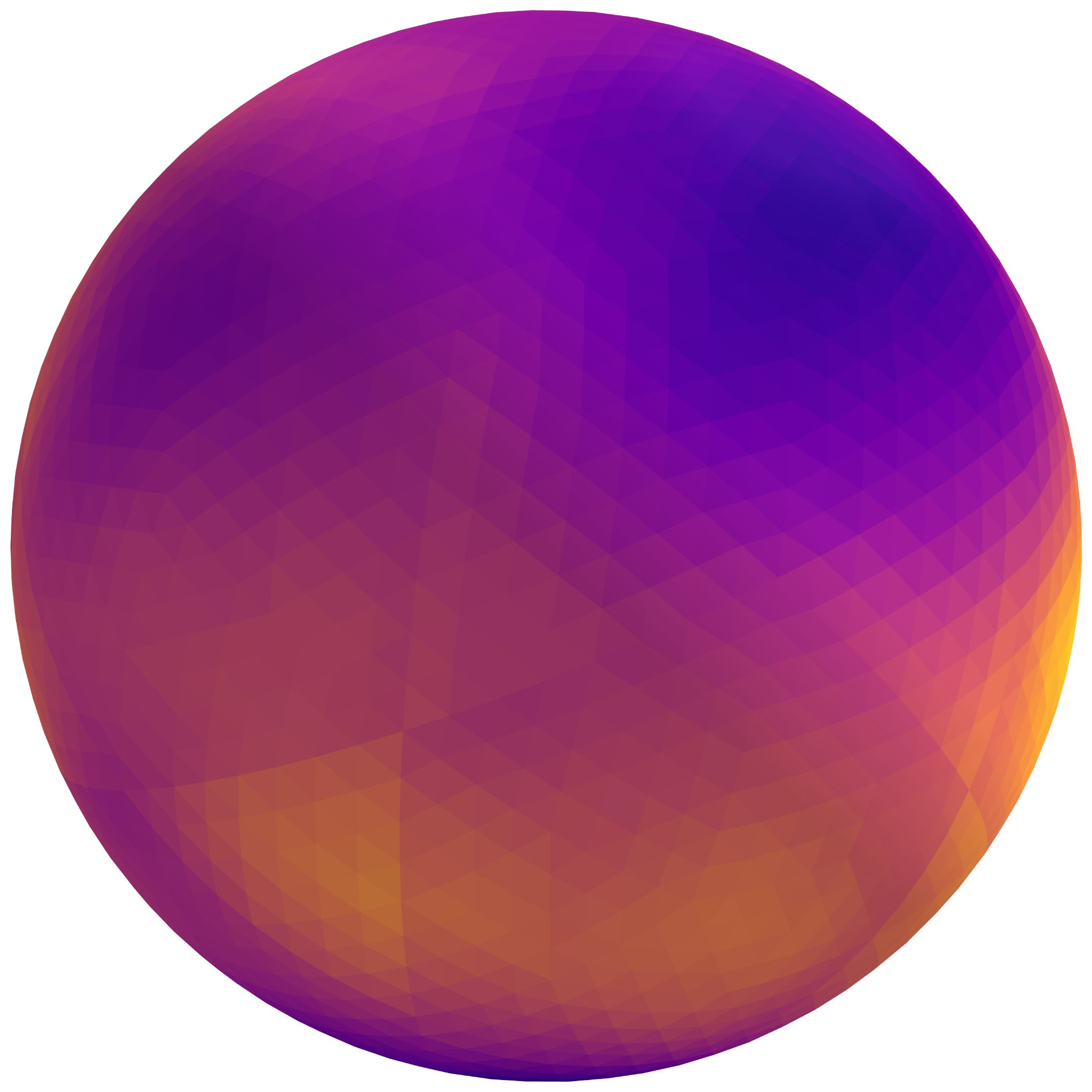}
            \caption{{\small Sample of $2$-cochain (function over faces)}}  
            \label{fig:2-sample}
        \end{subfigure}
        \vskip\baselineskip
        \begin{subfigure}[b]{\subfigwidth\textwidth}   
            \centering 
            \includegraphics[width=\textwidth]{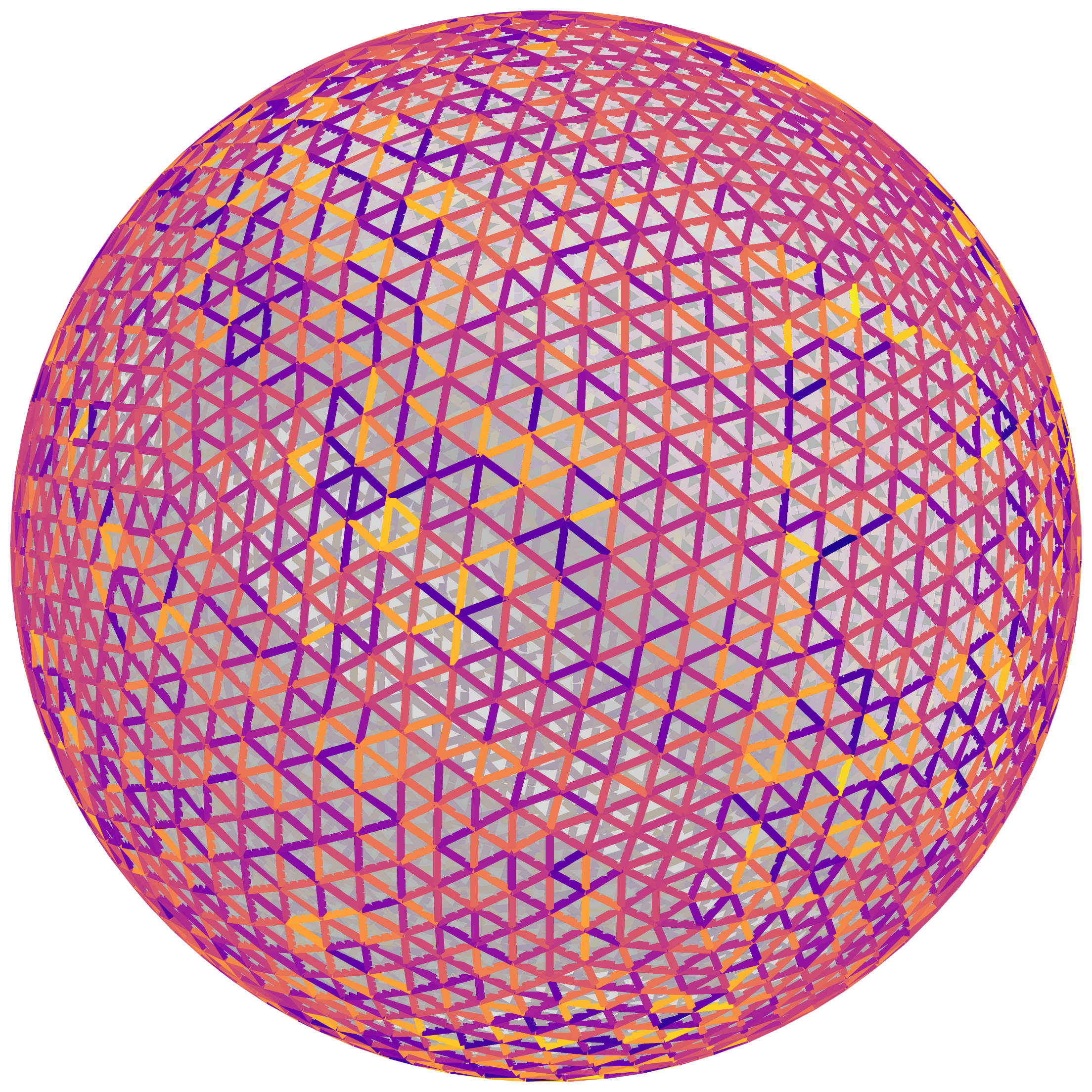}
            \caption[]{{\small Sample of $1$-cochain (function over edges) }}  
            \label{fig:1-sample}
        \end{subfigure}
        \hfill
        \begin{subfigure}[b]{\subfigwidth\textwidth}   
            \centering 
            \includegraphics[width=\textwidth]{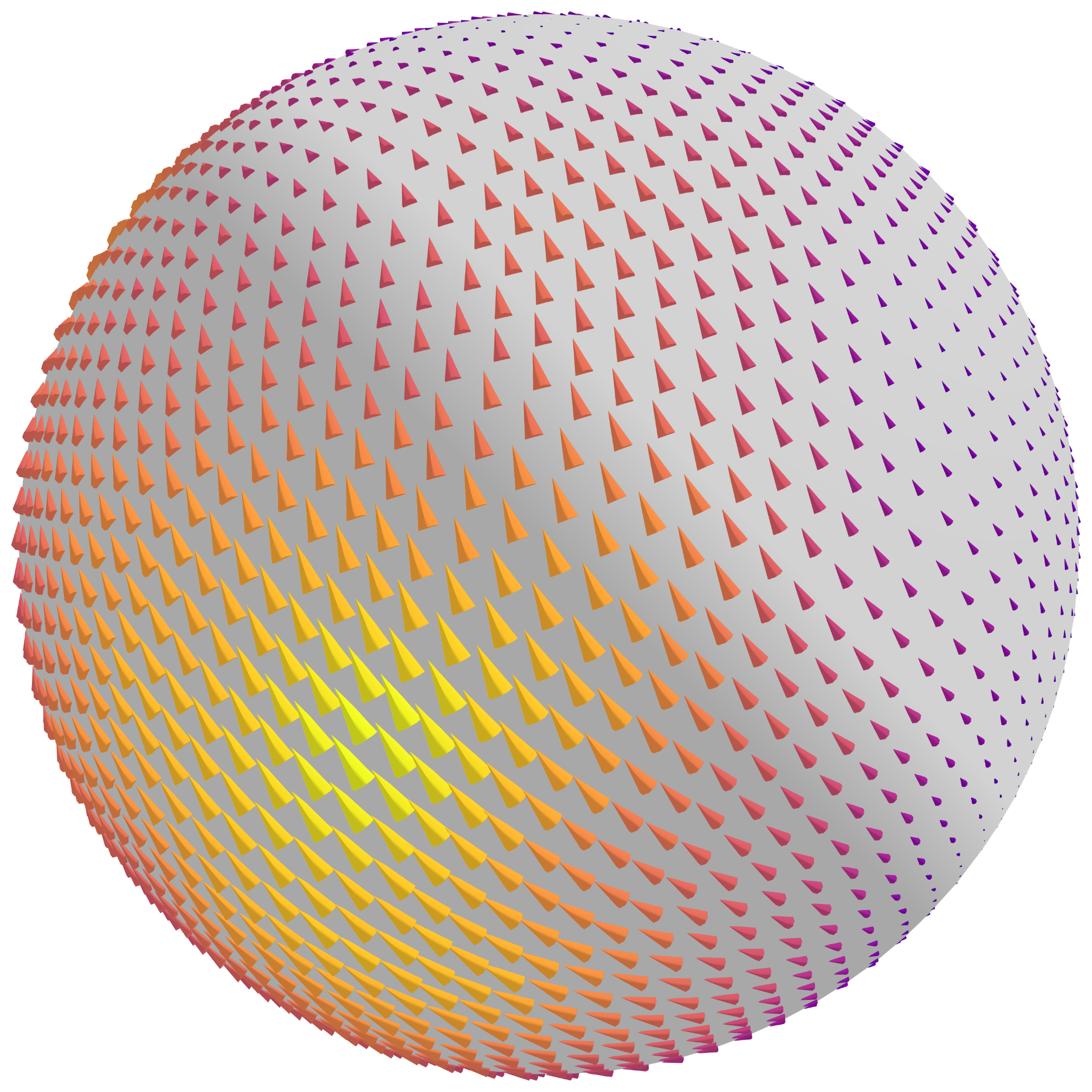}
            \caption{{\small Sample from Figure \ref{fig:1-sample}, projected into the tangent vector space (see Section \ref{sec:obs})}}  
            \label{fig:1-sample-vec}
        \end{subfigure}
        \caption{\small Samples of Whittle-Matérn fields on the $k$-cochains of a unit icosphere with 4 subdivisions. $\kappa = 2$ and $\alpha = 3$.} 
        \label{fig:samples}
    \end{figure*}

\subsection{Posterior Conditioning}\label{sec:post}

\begin{figure*}[h!]
        \centering
        \begin{subfigure}[b]{0.24\textwidth}
            \centering
            \includegraphics[width=\textwidth]{face_sample.png}
            \caption{{\small Sample of $2$-cochain\\\phantom{.}}}  
            \label{fig:2-sample-post}
        \end{subfigure}
        \hfill
        \begin{subfigure}[b]{0.24\textwidth}  
            \centering 
            \includegraphics[width=\textwidth]{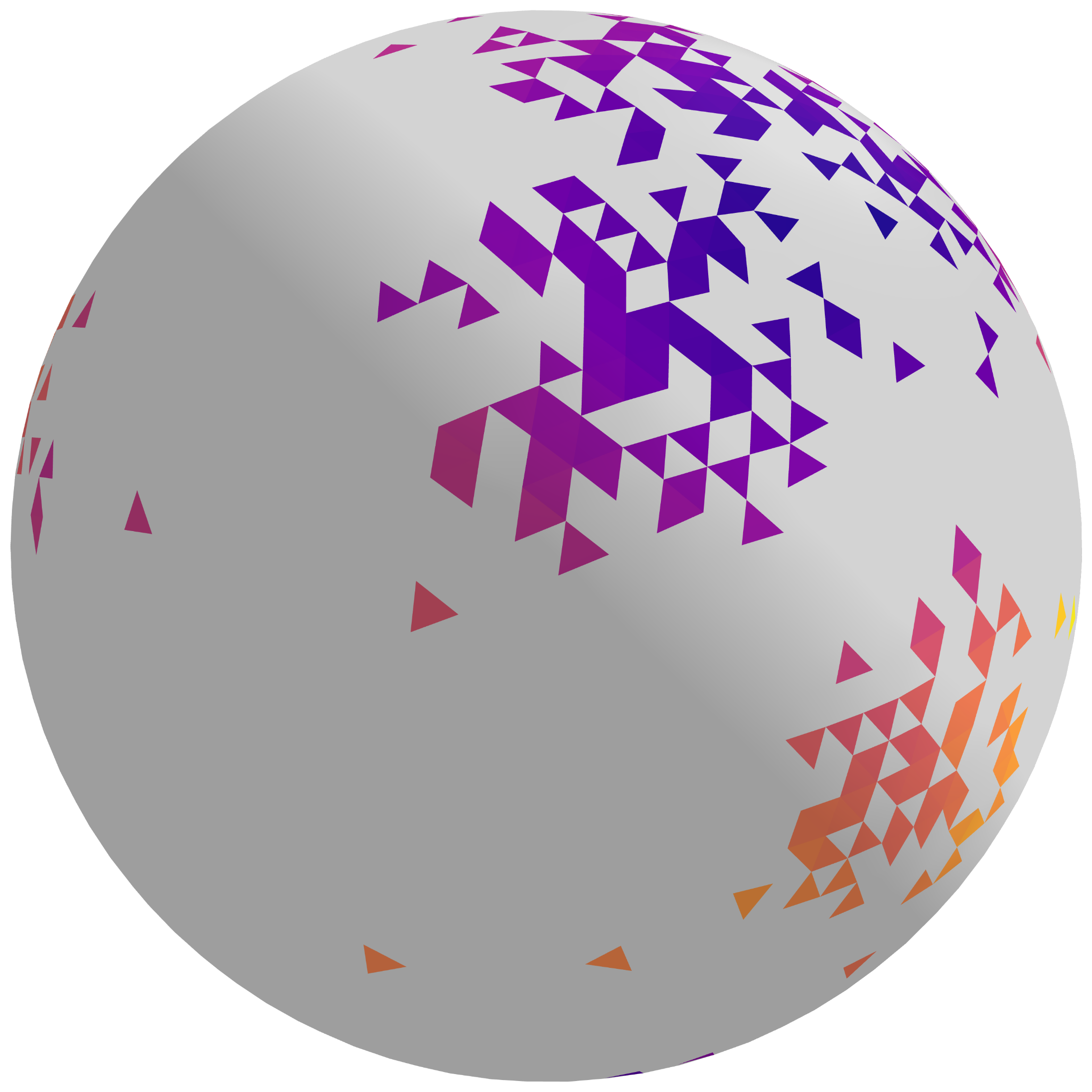}
            \caption{{\small 500 observed points\\\phantom{.}}}  
            \label{fig:obs-pts}
        \end{subfigure}
        \hfill
        \begin{subfigure}[b]{0.24\textwidth}   
            \centering 
            \includegraphics[width=\textwidth]{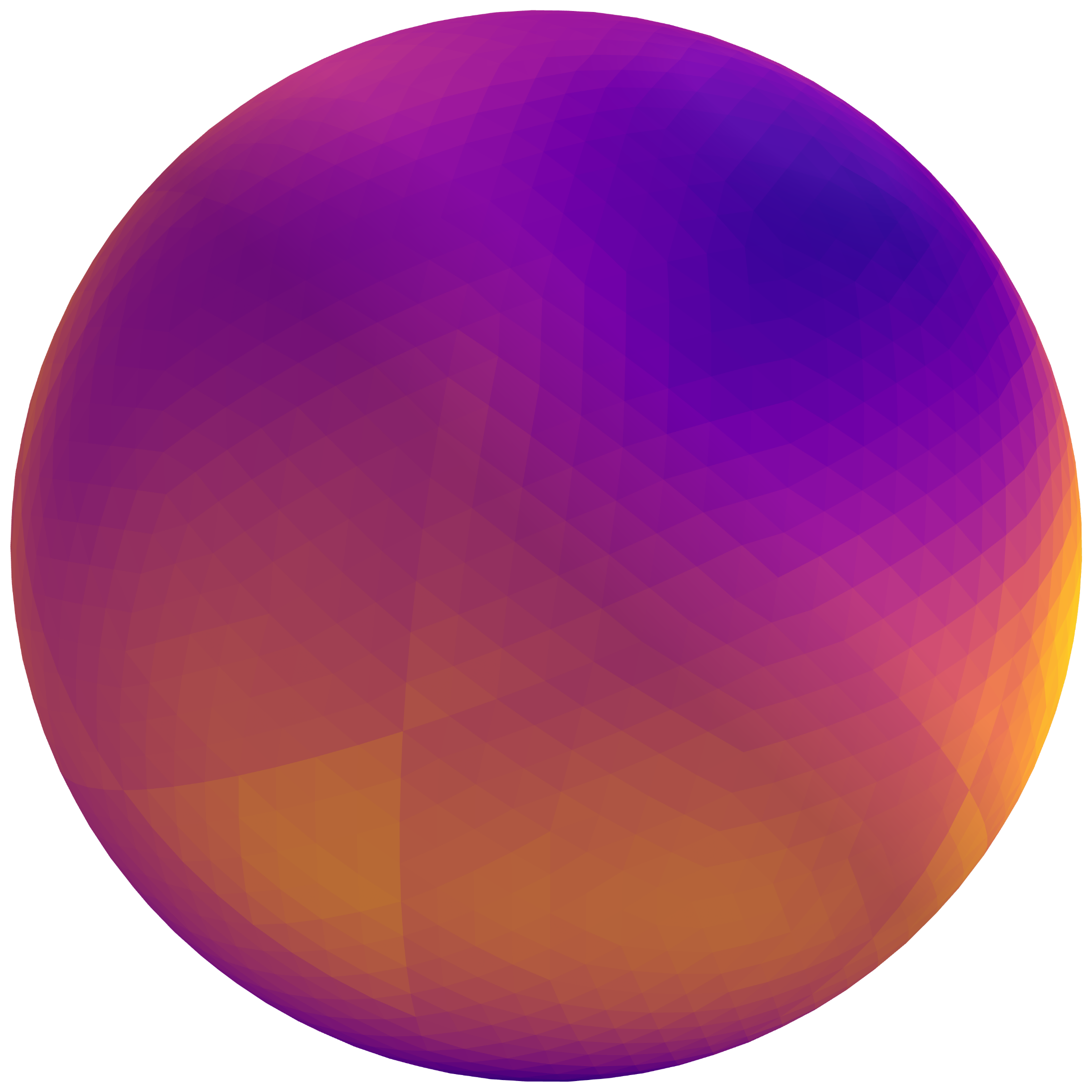}
            \caption[]{{\small Posterior mean\\\phantom{.}}}  
            \label{fig:post-mean}
        \end{subfigure}
        \hfill
        \begin{subfigure}[b]{0.24\textwidth}   
            \centering 
            \includegraphics[height=\textwidth]{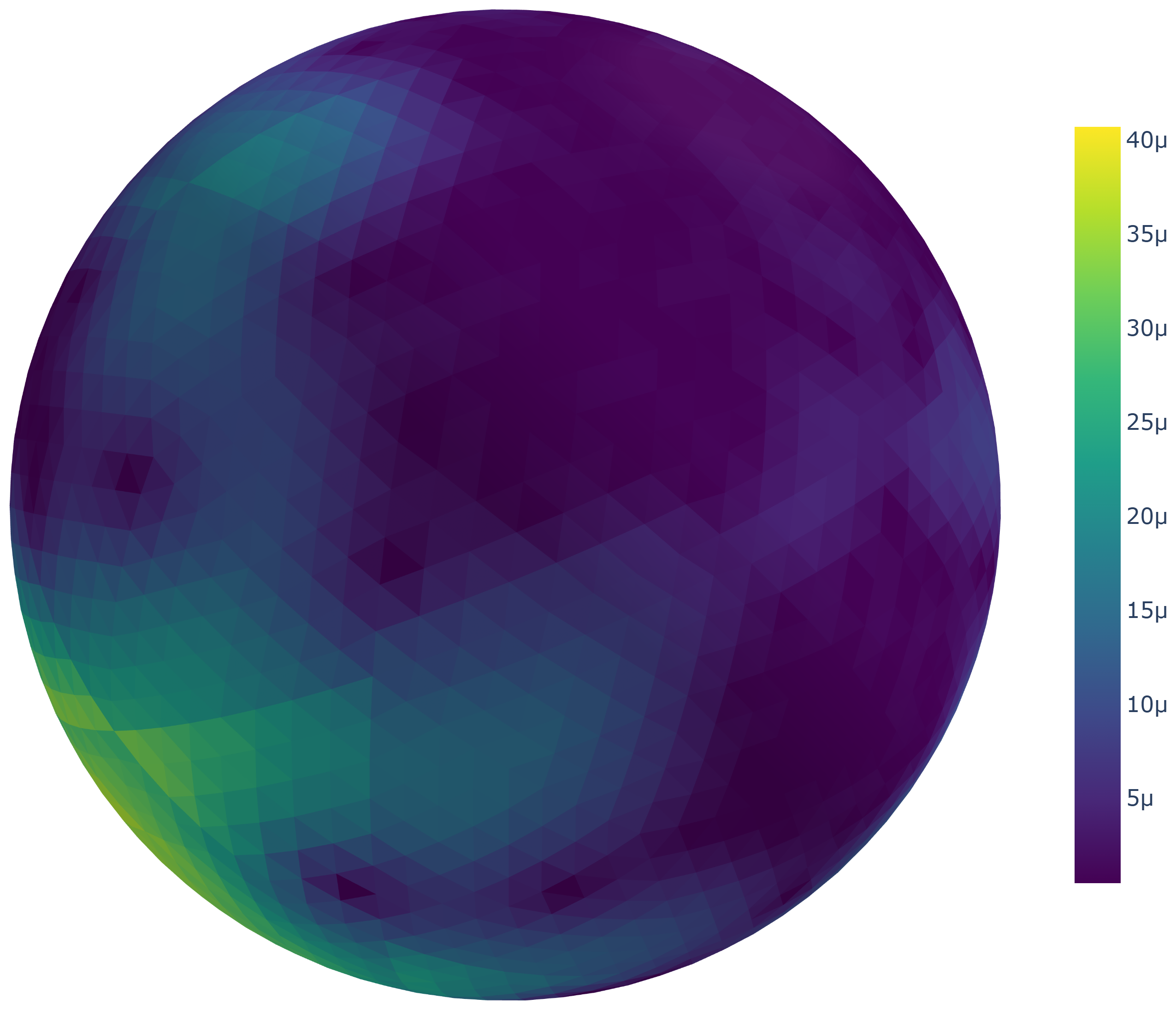}
            \caption{{\small Posterior standard deviation}}  
            \label{fig:post-std}
        \end{subfigure}
        \caption{\small Posterior conditioning of the Whittle-Matérn field on a $2$-cochain. RMSE = $1.24 \times 10^{-5}$} 
        \label{fig:posterior}
    \end{figure*}

The posterior covariance $\text{Cov}(x \mid y) = Q_{\text{post}}^{-1}$ where 
\begin{equation}
    Q_\text{post} = Q_\alpha + A^\top R^{-1} A.
\end{equation}
$A$ is an observation operator that maps the latent cochain to the observations, $y = Ax + \epsilon$, and $R$ is the covariance of the observation noise, $\epsilon \sim \mathcal{N}(0, R)$ e.g. for i.i.d noise with variance $\tau^2$, $R = \tau^2I$. For $R = \tau^2I$, the posterior mean of $x$ is given by
\begin{equation}
\mathbb{E}[x \mid y] = Q_\text{post}^{-1} \left(\frac{A^\top y}{\tau^2} \right).
\end{equation}
$Q_\text{post}^{-1}$ is a dense matrix, but we can find the mean and take posterior samples without ever forming it explicitly. For the mean we employ two triangular solves
\begin{equation}
    Q_\text{post} = L_\text{post}L_\text{post}^\top, \quad L_\text{post}w = A^\top R^{-1} y, \quad L_\text{post}^\top\mu = w.
\end{equation}
For posterior samples, we employ one additional triangular solve,
\begin{equation}
    x = \mu + L_\text{post}^{-\top} z, \quad z \sim \mathcal{N}(0, I).
\end{equation}
We can find the posterior marginal variance approximately through the Monte Carlo sample variance, or exactly by solving Takahashi equations to find only $\text{diag}(Q_\text{post}^{-1})$ \citep{rueGaussianMarkovRandom2005}. Figure \ref{fig:posterior} shows posterior inference on a $2$-cochain on a sphere.

\subsection{Observation Operators}\label{sec:obs}

The latent field is a cochain $x\in C^k(K)$ with prior $\mathcal{N}(0,Q_\alpha^{-1})$. Any observation that is a linear functional of the underlying form (e.g. an integral, a point value, or a vector) enters through the observation operator $A$ in $y=Ax+\epsilon$, and the posterior of Section~\ref{sec:post} applies unchanged.

\paragraph{Mesh-aligned observations.} If $y_i$ is an integral over a chain $c\in C_k(K)$, then by the de Rham map $\int_c\omega=\langle\psi^k\omega,c\rangle$, so the corresponding row of $A$ is the coefficient vector of $c$. For a single simplex this is a signed indicator, with sign given by the relative orientation. Note that this includes scalar values observed at vertices of the mesh.

\paragraph{Point observations.} For non-mesh aligned observations, we evaluate the Whitney map $\phi^kx$ at the observation location $p$ as in Eq. \ref{eq:whitney-map}, using the barycentric coordinates $\lambda_v(p)$ of $p$. Only the Whitney forms of the face $f$ containing $p$ contribute, so each row of $A$ is sparse. For $k=1$, when observations are vectors rather than 1-forms, we use a vector proxy for $\phi^1 x$,
\begin{equation}\label{eq:whitney-obs}
  \mathbf{v}(p)=\sum_{[v_a,v_b]\subset\partial f}x([v_a,v_b])
  \big(\lambda_a(p)\nabla\lambda_b-\lambda_b(p)\nabla\lambda_a\big),
\end{equation}
where $\nabla\lambda_i$ is constant on $f$. An observed
vector is compared with $E_p^\top\mathbf{v}(p)$, where
$E_p\in\mathbb{R}^{3\times2}$ is the observation's local frame.

\paragraph{Vector-field priors.} For a collection of points $p_1, \dots, p_m$, we can define the operator $B:C^1(K)\to\mathbb{R}^{2m}$ by evaluating $E_{p_i}^\top\mathbf{v}(p_i)$ at each $p_i$. The prior on vector fields is the push-forward of the cochain prior through $B$, $\mathbf{v}=Bx\sim\mathcal{N}(0,BQ_\alpha^{-1}B^\top)$. We can avoid constructing this dense covariance explicitly; samples are obtained by sampling $x$ as in Section \ref{sec:prior} and applying $B$, and vector observations enter inference through the corresponding rows of $B$ in $A$.

\section{Experiments}

As a demonstration, we apply this model to interpolation of a vector field on a sphere. The data are monthly averaged $u$ and $v$ components of wind at the 500hPa pressure level for the month of July 2020, taken from the ERA5 atmospheric reanalysis dataset \citep{hersbachERA5MonthlyAveraged2023}\footnote{Strictly speaking, these model outputs are the spatial averages over $0.25^\circ$ grid cells, but here we treat them as point estimates.}, scaled to have an average norm of 1. For observations we subsample the data to a $10^\circ \times 10^\circ$ grid, and infer back to the $0.25^\circ \times 0.25^\circ$ grid, giving 1,035,360 output locations. The ground truth data is shown in Figure \ref{fig:ground-truth-vecs} and the observations are shown in Figure \ref{fig:obs-vecs}.

\newcommand{\figwidth}{0.65}

\begin{figure}
    \centering
    \includegraphics[width=\figwidth\linewidth]{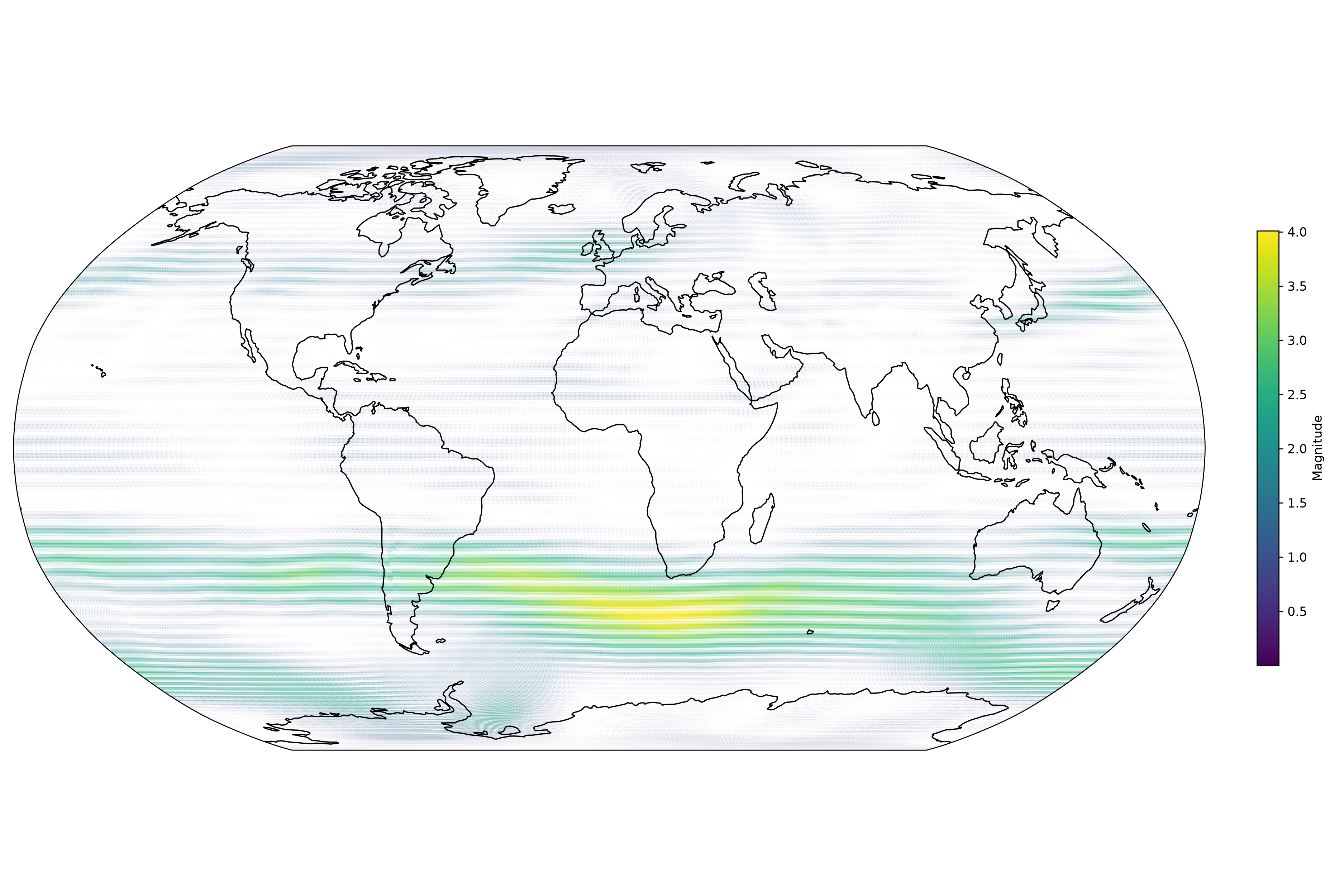}
    \caption{Monthly average wind data from ERA-5 on a $0.25^\circ \times 0.25^\circ$ grid, shown on the Robinson projection. Arrow colour and size indicates magnitude.}
    \label{fig:ground-truth-vecs}
\end{figure}
\begin{figure}
    \centering
    \includegraphics[width=\figwidth\linewidth]{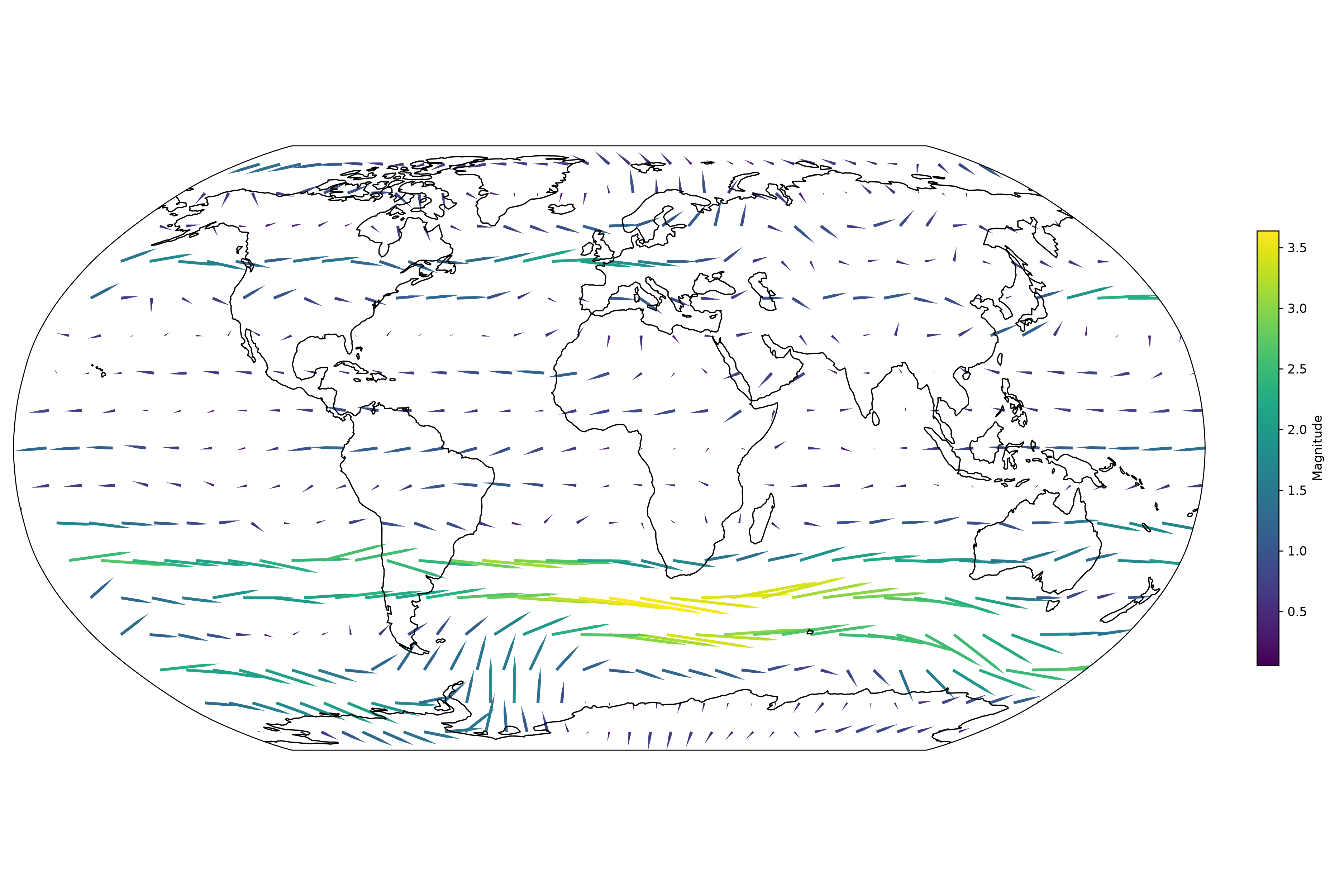}
    \caption{Observed wind vectors, downsampled from ground truth to a $10^\circ \times 10^\circ$ grid, shown on the Robinson projection. Arrow colour and size indicates magnitude.}
    \label{fig:obs-vecs}
\end{figure}

We begin by generating an icosphere mesh with 5 subdivisions, producing a well-centred mesh with 10242 vertices, 30720 edges, and 20480 faces, and using PyDEC we assemble a simplicial complex. For hyperparameter values, we fix $\alpha = 2 \ (\nu = 1)$ and estimate $\rho$ by empirically measuring the spatial correlation of the observations. For our data this gives $\rho = 0.615\operatorname{rad}$ and so $\kappa = 4.597$.

\begin{figure}[h!]
    \centering
    \includegraphics[width=0.5\linewidth]{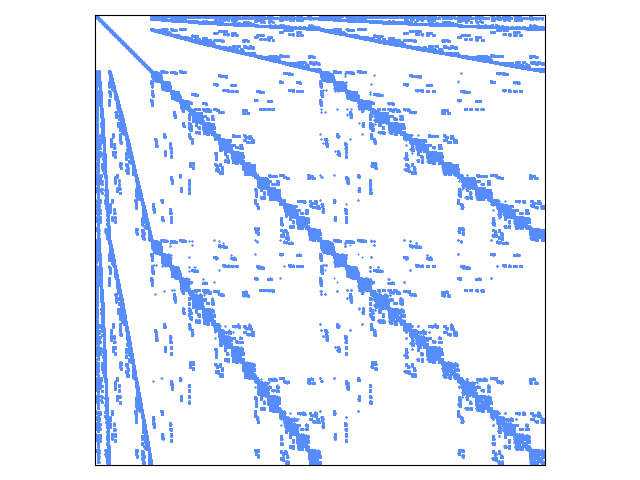}
    \caption{Sparsity pattern of $Q_2$ when $k=1$. Matrix dimensions are 30720 $\times$ 30720.}
    \label{fig:mat-sparse}
\end{figure}

Using the DEC operators defined by the simplicial complex, we construct $Q_2$ for $k=1$ as in Eq. \ref{eq:q_2}, giving a prior precision over the edges of the mesh. The sparsity pattern of this matrix is shown in Figure \ref{fig:mat-sparse}. We then construct our observation operator $A$ as the Whitney map for vector field priors described in Section \ref{sec:obs}. Note that since we use an icosphere, neither the observations nor the ground truth are (necessarily) aligned with the vertices of the mesh. We set $\tau = 0.001$ and use $R= \tau^2 I$ for the observation noise.

\begin{figure}[h!]
    \centering
    \includegraphics[width=\figwidth\linewidth]{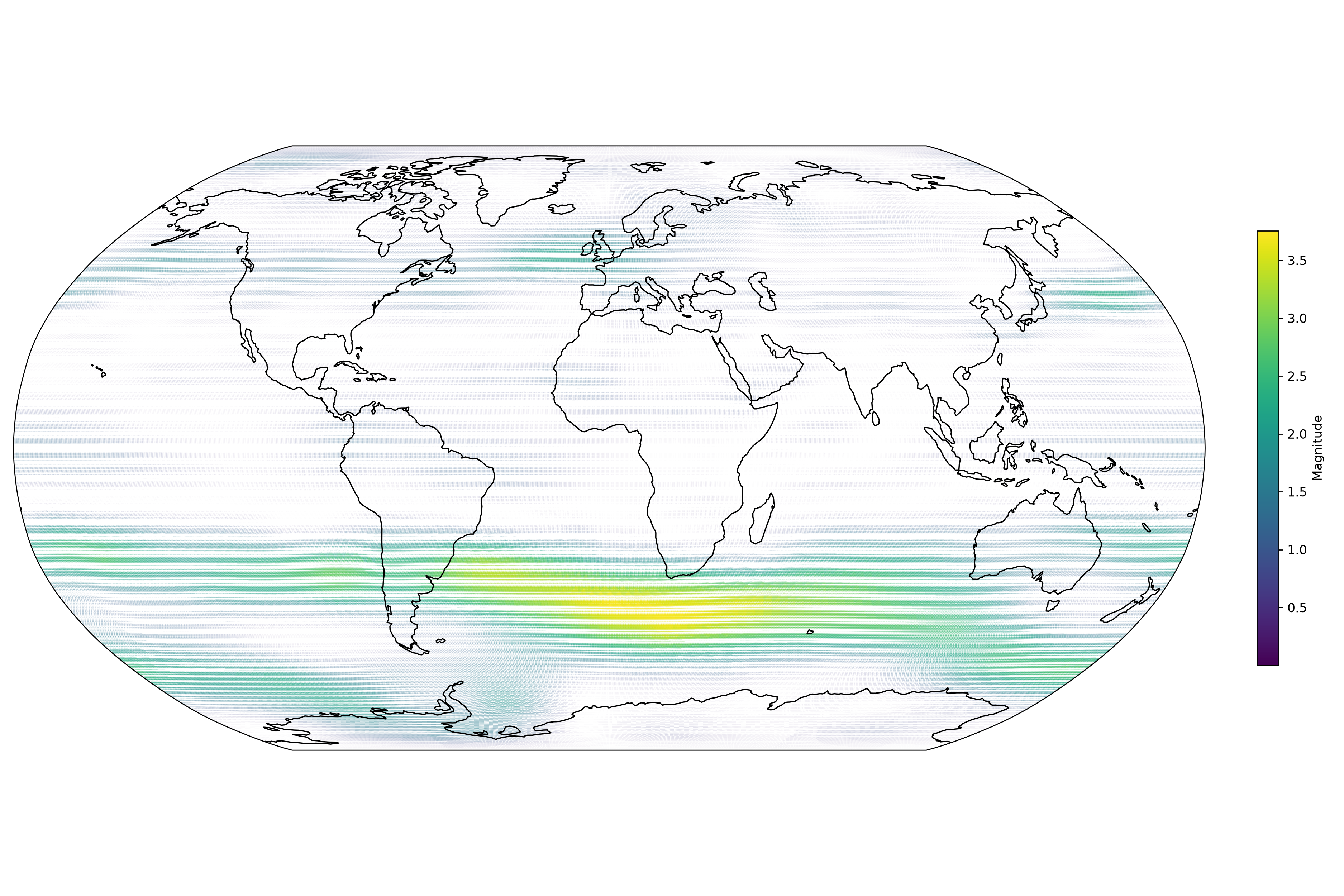}
    \caption{Posterior mean field, $0.25^\circ \times 0.25^\circ$, conditioned on $10^\circ \times 10^\circ$ grid of observation vectors. Arrow colour and size indicates magnitude. MSE = 0.046}
    \label{fig:post-mean-vecs}
\end{figure}

\begin{figure}[h!]
    \centering
    \includegraphics[width=\figwidth\linewidth]{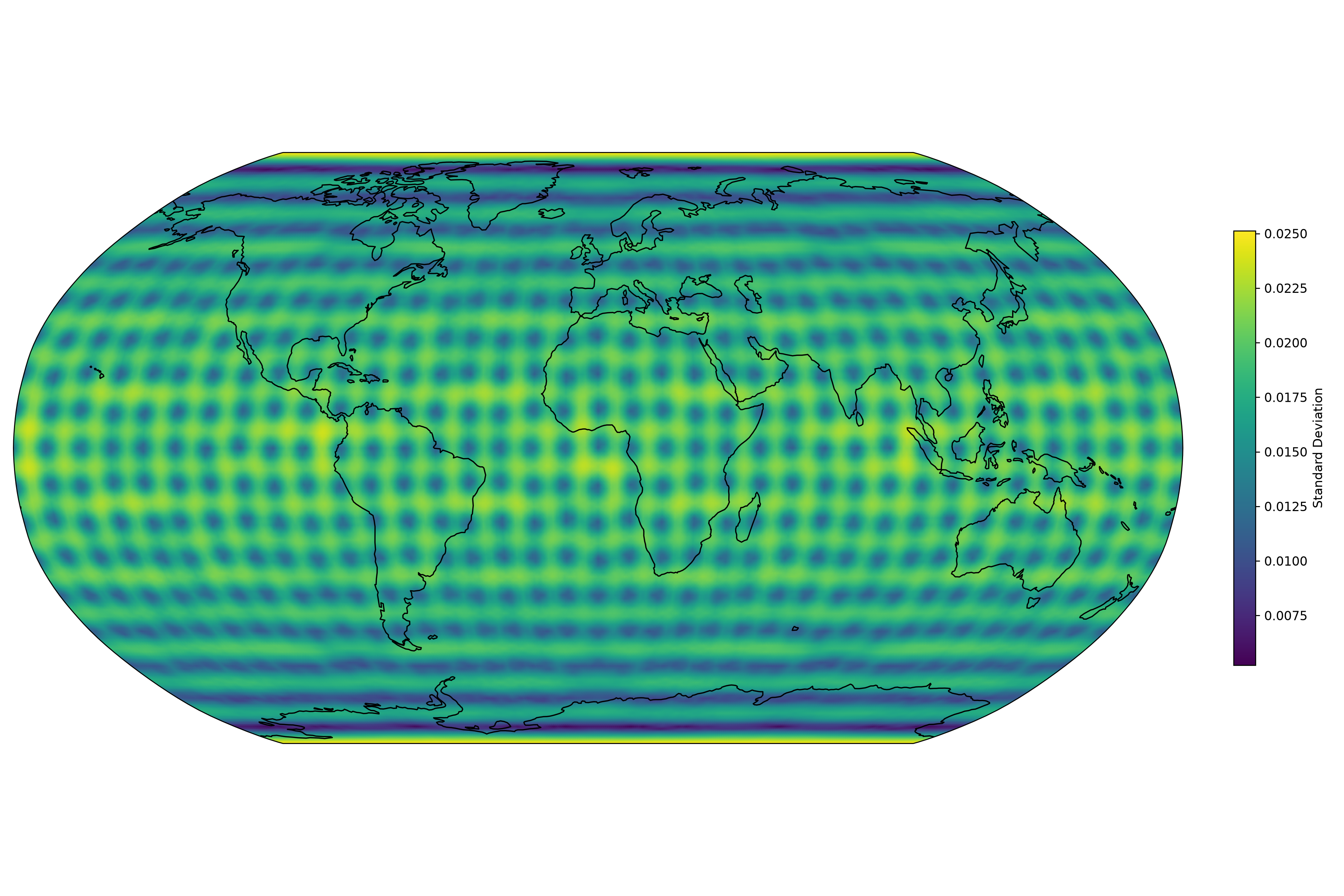}
    \caption{Posterior standard deviation of wind vector \emph{along} the edges of the mesh, interpolated to lat/lon grid using $\phi^0$, the Whitney map for when $k=0$. Colour indicates size of standard deviation.}
    \label{fig:post-std-vecs}
\end{figure}

\begin{figure}[h!]
    \centering
    \includegraphics[width=\figwidth\linewidth]{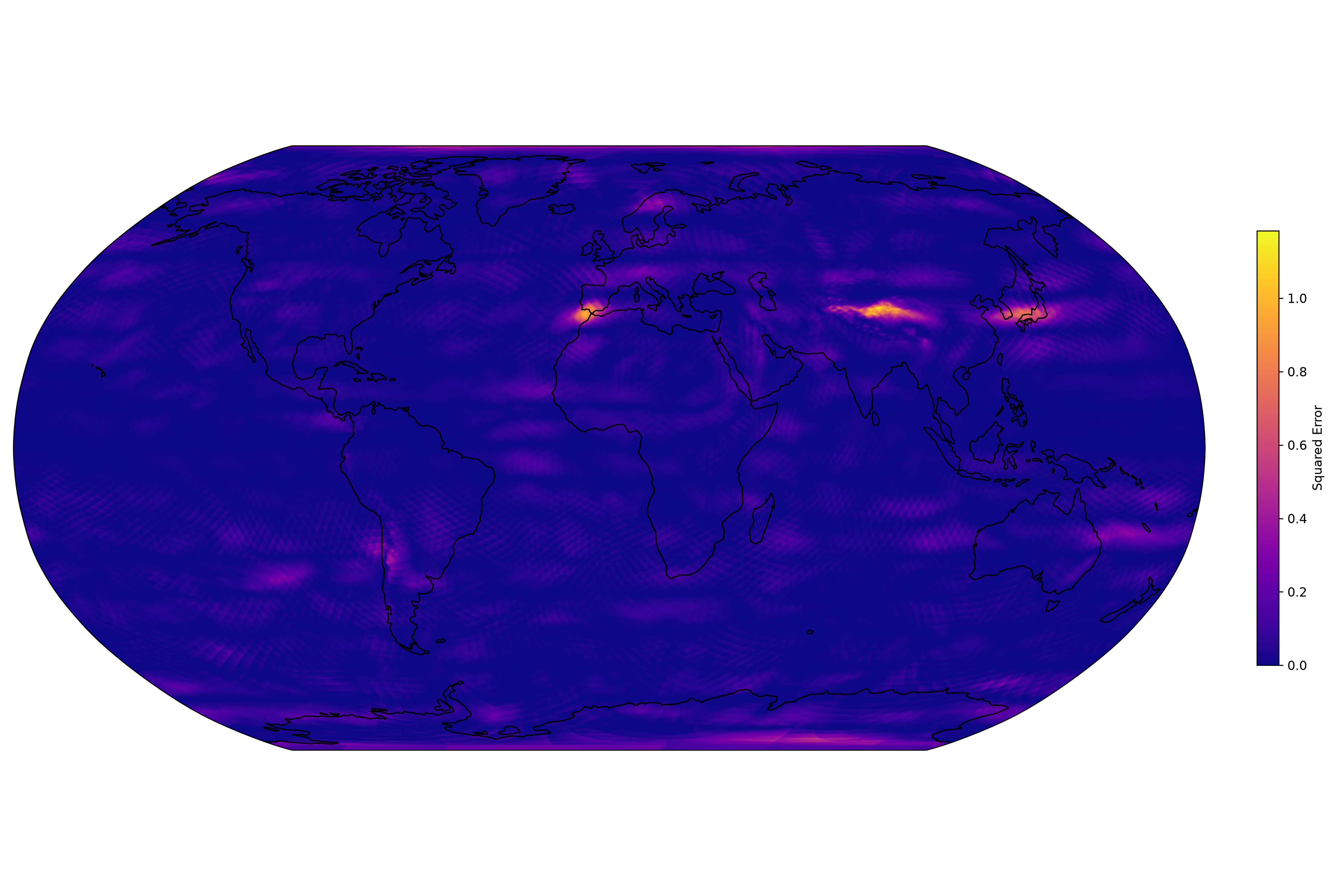}
    \caption{Square error of posterior mean field with respect to ground truth. Colour indicates magnitude of error. Squared error is calculated as $||\mathbf{x} - \mu_{\mathbf{x|y}}||^2_2$.}
    \label{fig:post-error-vecs}
\end{figure}

Using the procedure detailed in Section \ref{sec:post}, we obtain the posterior mean and standard deviation, shown in Figures \ref{fig:post-mean-vecs} and \ref{fig:post-std-vecs}, respectively. The posterior estimate has a mean square error of 0.046 with respect to the ground truth, and the spatial distribution of square error is shown in Figure \ref{fig:post-error-vecs}.


\section{Conclusion}

This work presents a sparse precision formulation for Whittle-Matérn fields for $k$-cochains using DEC operators. By applying Whitney maps and introducing the Whitney observation operators, this model can be applied to any point, edge, areal, or vector field data, with observation and prediction locations anywhere on the underlying manifold. Future development will focus on boundary conditions and non-stationarity, as well as sparse approximations of Hodge-compositional covariance kernels.

\clearpage

\bibliographystyle{chicago}
\bibliography{MeshGMRF}

@inproceedings{alain2026probabilistic,
  title = {Probabilistic Modelling of Vector Fields Using Structure-Preserving Discretisation},
  booktitle = {{{NeurIPS}} 2026 Workshop on {{AI}} for Stochastic Dynamics: {{From}} Theoretical Foundations to Scientific Applications},
  author = {Alain, Mathieu and Takao, So},
  year = 2026
}

@misc{alainGaussianProcessesCellular2024,
  title = {Gaussian {{Processes}} on {{Cellular Complexes}}},
  author = {Alain, Mathieu and Takao, So and Paige, Brooks and Deisenroth, Marc Peter},
  year = 2024,
  month = aug,
  number = {arXiv:2311.01198},
  eprint = {2311.01198},
  primaryclass = {cs.LG},
  publisher = {arXiv},
  doi = {10.48550/arXiv.2311.01198},
  urldate = {2026-10-20},
  archiveprefix = {arXiv}
}

@article{bellPyDECSoftwareAlgorithms2012,
  title = {{{PyDEC}}: {{Software}} and {{Algorithms}} for {{Discretization}} of {{Exterior Calculus}}},
  shorttitle = {{{PyDEC}}},
  author = {Bell, Nathan and Hirani, Anil N.},
  year = 2012,
  month = nov,
  journal = {ACM Trans. Math. Softw.},
  volume = {39},
  number = {1},
  pages = {3:1--3:41},
  issn = {0098-3500},
  doi = {10.1145/2382585.2382588},
  urldate = {2024-10-24}
}

@article{bolinRationalSPDEApproach2020,
  title = {The {{Rational SPDE Approach}} for {{Gaussian Random Fields With General Smoothness}}},
  author = {Bolin, David and Kirchner, Kristin},
  year = 2020,
  month = apr,
  journal = {Journal of Computational and Graphical Statistics},
  volume = {29},
  number = {2},
  pages = {274--285},
  issn = {1061-8600, 1537-2715},
  doi = {10.1080/10618600.2019.1665537},
  urldate = {2026-10-08},
  langid = {english}
}

@inproceedings{borovitskiyMaternGaussianProcesses2020,
  title = {Mat\'ern {{Gaussian Processes}} on {{Riemannian Manifolds}}},
  booktitle = {Advances in {{Neural Information Processing Systems}}},
  author = {Borovitskiy, Viacheslav and Terenin, Alexander and Mostowsky, Peter and Deisenroth (he/him), Marc},
  year = 2020,
  volume = {33},
  pages = {12426--12437},
  publisher = {Curran Associates, Inc.},
  urldate = {2026-06-10}
}

@inproceedings{borovitskiyMaternGaussianProcesses2021b,
  title = {Mat\'ern {{Gaussian Processes}} on {{Graphs}}},
  booktitle = {Proceedings of {{The}} 24th {{International Conference}} on {{Artificial Intelligence}} and {{Statistics}}},
  author = {Borovitskiy, Viacheslav and Azangulov, Iskander and Terenin, Alexander and Mostowsky, Peter and Deisenroth, Marc and Durrande, Nicolas},
  year = 2021,
  month = mar,
  pages = {2593--2601},
  publisher = {PMLR},
  issn = {2640-3498},
  urldate = {2026-10-20},
  langid = {english}
}

@misc{desbrunDiscreteExteriorCalculus2005,
  title = {Discrete {{Exterior Calculus}}},
  author = {Desbrun, Mathieu and Hirani, Anil N. and Leok, Melvin and Marsden, Jerrold E.},
  year = 2005,
  month = aug,
  number = {arXiv:math/0508341},
  eprint = {math/0508341},
  publisher = {arXiv},
  doi = {10.48550/arXiv.math/0508341},
  urldate = {2026-09-25},
  archiveprefix = {arXiv}
}

@book{flemingFunctionsSeveralVariables1987,
  title = {Functions of {{Several Variables}}},
  author = {Fleming, Wendell H.},
  year = 1987,
  month = jun,
  publisher = {Springer Science \& Business Media},
  googlebooks = {xK0Crp3kt8EC},
  isbn = {978-0-387-90206-7},
  langid = {english}
}

@article{gillanDiscreteGaussianVector2026,
  title = {Discrete {{Gaussian}} Vector Fields on Meshes},
  author = {Gillan, Michael and Siegert, Stefan and Youngman, Benjamin D.},
  year = 2026,
  month = jun,
  journal = {Stochastic Environmental Research and Risk Assessment},
  volume = {40},
  number = {6},
  pages = {161},
  issn = {1436-3259},
  doi = {10.1007/s00477-026-03292-1}
}

@misc{hersbachERA5MonthlyAveraged2023,
  title = {{{ERA5}} Monthly Averaged Data on Pressure Levels from 1940 to Present},
  author = {Hersbach, H. and Bell, B. and Berrisford, P. and Biavati, G. and Hor{\'a}nyi, A. and Mu{\~n}oz Sabater, J. and Nicolas, J. and Peubey, C. and Radu, R. and Rozum, I. and Schepers, D. and Simmons, A. and Soci, C. and Dee, D. and Th{\'e}paut, J-N},
  year = 2023,
  publisher = {Copernicus Climate Change Service (C3S) Climate Data Store (CDS)},
  doi = {10.24381/CDS.6860A573},
  urldate = {2025-05-23}
}

@phdthesis{hiraniDiscreteExteriorCalculus2003a,
  title = {Discrete Exterior Calculus},
  author = {Hirani, Anil Nirmal},
  year = 2003,
  urldate = {2026-06-10},
  school = {California Institute of Technology}
}

@inproceedings{hutchinsonVectorvaluedGaussianProcesses2021a,
  title = {Vector-Valued {{Gaussian}} Processes on Riemannian Manifolds via Gauge Independant Projected Kernels},
  booktitle = {Proceedings of the 35th {{International Conference}} on {{Neural Information Processing Systems}}},
  author = {Hutchinson, Michael and Terenin, Alexander and Borovitskiy, Viacheslav and Takao, So and Teh, Yee Whye and Deisenroth, Marc Peter},
  year = 2021,
  month = dec,
  series = {{{NIPS}} '21},
  pages = {17160--17169},
  publisher = {Curran Associates Inc.},
  address = {Red Hook, NY, USA},
  urldate = {2026-06-10},
  isbn = {978-1-7138-4539-3}
}

@article{lindgrenExplicitLinkGaussian2011,
  title = {An Explicit Link between {{Gaussian}} Fields and {{Gaussian Markov}} Random Fields: The Stochastic Partial Differential Equation Approach},
  shorttitle = {An Explicit Link between {{Gaussian}} Fields and {{Gaussian Markov}} Random Fields},
  author = {Lindgren, Finn and Rue, H{\aa}vard and Lindstr{\"o}m, Johan},
  year = 2011,
  journal = {Journal of the Royal Statistical Society: Series B (Statistical Methodology)},
  volume = {73},
  number = {4},
  pages = {423--498},
  issn = {1467-9868},
  doi = {10.1111/j.1467-9868.2011.00777.x},
  urldate = {2024-02-19},
  langid = {english}
}

@article{lindgrenSPDEApproachGaussian2022,
  title = {The {{SPDE}} Approach for {{Gaussian}} and Non-{{Gaussian}} Fields: 10 Years and Still Running},
  shorttitle = {The {{SPDE}} Approach for {{Gaussian}} and Non-{{Gaussian}} Fields},
  author = {Lindgren, Finn and Bolin, David and Rue, H{\aa}vard},
  year = 2022,
  month = aug,
  journal = {Spatial Statistics},
  volume = {50},
  eprint = {2111.01084},
  primaryclass = {stat},
  pages = {100599},
  issn = {22116753},
  doi = {10.1016/j.spasta.2022.100599},
  urldate = {2024-02-19},
  archiveprefix = {arXiv}
}

@misc{nambirajanApproximatingGaussianWhittleMatern2026,
  title = {Approximating {{Gaussian Whittle-Matern Fields}} over {{Well-Centered Triangulations}} of {{Riemannian Manifolds}}},
  author = {Nambirajan, Srinivas},
  year = 2026,
  month = jun,
  number = {arXiv:2606.13827},
  eprint = {2606.13827},
  primaryclass = {math.NA},
  publisher = {arXiv},
  doi = {10.48550/arXiv.2606.13827},
  urldate = {2026-08-21},
  archiveprefix = {arXiv}
}

@book{rasmussenGaussianProcessesMachine2006,
  title = {Gaussian Processes for Machine Learning},
  author = {Rasmussen, Carl Edward and Williams, Christopher K. I.},
  year = 2006,
  series = {Adaptive Computation and Machine Learning},
  publisher = {MIT Press},
  address = {Cambridge, Mass},
  isbn = {978-0-262-18253-9},
  lccn = {QA274.4 .R37 2006}
}

@inproceedings{robert-nicoudIntrinsicGaussianVector2024b,
  title = {Intrinsic {{Gaussian Vector Fields}} on {{Manifolds}}},
  booktitle = {Proceedings of {{The}} 27th {{International Conference}} on {{Artificial Intelligence}} and {{Statistics}}},
  author = {{Robert-Nicoud}, Daniel and Krause, Andreas and Borovitskiy, Viacheslav},
  year = 2024,
  month = apr,
  pages = {1306--1314},
  publisher = {PMLR},
  issn = {2640-3498},
  urldate = {2026-06-10},
  langid = {english}
}

@book{rueGaussianMarkovRandom2005,
  title = {Gaussian {{Markov Random Fields}}: {{Theory}} and {{Applications}}},
  author = {Rue, Havard and Held, Leonhard},
  year = 2005,
  series = {Monographs on Statistics and Applied Probability},
  edition = {1},
  volume = {104},
  publisher = {{Chapman and Hall/CRC}},
  address = {Milton},
  doi = {10.1201/9780203492024},
  isbn = {978-1-58488-432-3}
}

@book{whitneyGeometricIntegrationTheory1957,
  title = {Geometric {{Integration Theory}}},
  author = {Whitney, Hassler},
  year = 1957,
  eprint = {j.ctt183q1bm},
  eprinttype = {jstor},
  publisher = {Princeton University Press},
  urldate = {2026-10-20},
  isbn = {978-0-691-07972-1}
}

@article{whittleStochasticProcessesSeveral1963,
  title = {Stochastic {{Processes}} in {{Several Dimensions}}},
  author = {Whittle, P},
  year = 1963,
  journal = {Bulletin of the International Statistical Institute},
  volume = {40},
  number = {2},
  pages = {974--994}
}

@inproceedings{yangHodgeCompositionalEdgeGaussian2024a,
  title = {Hodge-{{Compositional Edge Gaussian Processes}}},
  booktitle = {Proceedings of {{The}} 27th {{International Conference}} on {{Artificial Intelligence}} and {{Statistics}}},
  author = {Yang, Maosheng and Borovitskiy, Viacheslav and Isufi, Elvin},
  year = 2024,
  month = apr,
  pages = {3754--3762},
  publisher = {PMLR},
  issn = {2640-3498},
  urldate = {2026-06-10},
  langid = {english}
}

\end{document}